\documentclass{elsart3-1}

\usepackage{amssymb}

\usepackage[english]{babel}
\usepackage{amsmath,amsfonts,amssymb}
\usepackage{color}
\usepackage{mathabx}
\usepackage{dsfont}

\newtheorem{theorem}{Theorem}[section]
\newtheorem{lemma}[theorem]{Lemma}
\newtheorem{e-proposition}[theorem]{Proposition}

\newtheorem{e-definition}[theorem]{Definition\rm}

\renewcommand{\theequation}{\arabic{equation}}
\def\og{\leavevmode\raise.3ex\hbox{$\scriptscriptstyle\langle\!\langle$~}}
\def\fg{\leavevmode\raise.3ex\hbox{~$\!\scriptscriptstyle\,\rangle\!\rangle$}}

\newcommand{\TT}{\mathbb{T}}
\newcommand{\RR}{\mathbb{R}}
\newcommand{\cF}{\mathcal{F}}

\newcommand{\cS}{\mathcal{S}}
\newcommand{\cA}{\mathcal{A}}
\newcommand{\dd}{\mathrm{d}}
\newcommand{\ve}{\varepsilon}
\newcommand{\eps}{\varepsilon}

\newcommand{\lla}{\left\langle}
\newcommand{\rra}{\right\rangle}

\newcommand{\ds}{\displaystyle}
\newcommand{\Nt}{\vvvert}

\journal{the Acad\'emie des sciences}
\begin{document}
    \centerline{}
    \begin{frontmatter}
        
        
        \selectlanguage{english}
        \title{A note on hypocoercivity for kinetic equations with heavy-tailed equilibrium}

        
        \selectlanguage{english}
        \author[NA]{Nathalie Ayi},
        \ead{nathalie.ayi@sorbonne-universite.fr}
        \author[MH]{Maxime Herda},
        \ead{maxime.herda@inria.fr}
        \author[HH]{H\'el\`ene Hivert} and 
        \ead{helene.hivert@ec-lyon.fr}
        \author[IT]{Isabelle Tristani}
        \ead{isabelle.tristani@ens.fr}
        
        \address[NA]{Sorbonne Universit\'e, Universit\'e de Paris, CNRS, Laboratoire Jacques-Louis Lions, 4 place Jussieu, 75005 Paris, France}
        \address[MH]{Inria, Univ. Lille, CNRS, UMR 8524 - Laboratoire Paul Painlev\'e, F-59000 Lille, France}
        \address[HH]{Univ. Lyon, \'Ecole centrale de Lyon, CNRS UMR 5208, Institut Camille Jordan, F-69134 \'Ecully, France}
        \address[IT]{DMA, \'Ecole Normale Sup\'erieure, CNRS, PSL Research University, 45 rue d'Ulm, 75005 Paris, France}
        
        
        \medskip
        \begin{center}
            {\small Received *****; accepted after revision +++++\\
                Presented by }
        \end{center}
        
        \begin{abstract}
            \selectlanguage{english}
            In this paper we are interested in the large time behavior of linear kinetic equations with heavy-tailed local equilibria. Our main contribution concerns the kinetic L\'evy-Fokker-Planck equation, for which we adapt hypocoercivity techniques in order to show that solutions converge exponentially fast to the global equilibrium. Compared to the classical kinetic Fokker-Planck equation, the issues here concern the lack of symmetry of the non-local L\'evy-Fokker-Planck operator and the understanding of its regularization properties. As a complementary related result, we also treat the case of the heavy-tailed BGK equation. 
            
            
            \vskip 0.5\baselineskip
            
            \noindent{\bf R\'esum\'e} \vskip 0.5\baselineskip \noindent
            {\bf Une note sur l'hypocoercivit\'e pour les \'equations cin\'etiques avec \'equilibres \`a queue lourde. }Dans cet article, on s'int\'eresse au comportement en temps long d'\'equations cin\'etiques lin\'eaires dont les \'equilibres locaux sont \`a queue lourde. Notre contribution principale concerne l'\'equation de L\'evy-Fokker-Planck cin\'etique, pour laquelle nous adaptons des techniques d'hypocoercivit\'e afin de d\'emontrer la convergence exponentielle des solutions vers un \'equilibre global. En comparant au cas de l'\'equation de Fokker-Planck cin\'etique classique, les enjeux ici sont li\'es au manque de sym\'etrie de l'op\'erateur non-local de L\'evy-Fokker-Planck et \`a la compr\'ehension de ses propri\'et\'es de r\'egularisation. En compl\'ement de notre analyse, nous traitons \'egalement le cas de l'\'equation de BGK \`a queue lourde. 
            
            
            
        \end{abstract}
    \end{frontmatter}
    %
    
    \selectlanguage{english}
    \section{Introduction}\label{s:introduction}
    
    We consider a distribution function $f\equiv f(t,x,v)$ depending on time $t\geq0$, position $x\in\TT^d$ and velocity~$v\in\RR^d$ which satisfies the fractional kinetic Fokker-Planck equation
    \begin{equation}\label{eq:kinfracFP}
        \partial_t f + v\cdot\nabla_x f\ =\ \nabla_v\cdot(vf)-(-\Delta_v)^{\alpha/2}f\,.
    \end{equation}
    Here we assume $\alpha\in(0,2)$ and the fractional Laplacian $-(-\Delta_v)^{\alpha/2}$ is such that for any Schwartz function~$g:\RR^d\to\RR$, one has $\cF((-\Delta_v)^{\alpha/2}g)(\xi)\ =\ |\xi|^{\alpha}\cF(g)(\xi)$ where $\cF(\cdot)$ denotes the Fourier transform. There are many equivalent definitions of the fractional Laplacian (see \cite{kwasnicki_2017_fractional}). Among them we shall use
    \begin{equation}\label{eq:fraclap_int1}
        (-\Delta_v)^{\alpha/2}g(v)\ =\ C_{d,\alpha}\,\mathrm{P.V.}\int_{\RR^d}\frac{g(v)-g(w)}{|v-w|^{d+\alpha}}\, \dd w\,,
    \end{equation}
    where $\mathrm{P.V.}$ stands for the principal value and the constant is given by $C_{d,\alpha}={2^\alpha\Gamma(\frac{d+\alpha}{2})}/{(\pi^{d/2}|\Gamma(-\frac{\alpha}{2})|)}$ where $\Gamma(\cdot)$ is the Gamma function. In the following we drop the principal value in the notations. We denote the L\'evy-Fokker-Planck operator appearing on the right-hand side of~\eqref{eq:kinfracFP} by 
    \[
    L_\alpha g\ =\ \nabla_v\cdot(v\,g)-(-\Delta_v)^{\alpha/2}g\,.
    \]
    By passing to Fourier variables one has $\cF(L_\alpha g)(\xi) = -\xi\cdot\nabla_\xi\hat{g}(\xi)-|\xi|^{\alpha}\hat{g}(\xi)$, where $\hat{g} = \cF(g)$. From this formula, one sees that the function 
    \[
    \mu_\alpha(v)\ =\  Z_{d,\alpha}^{-1}\cF^{-1}\left(e^{-|\xi|^\alpha /\alpha}\right)
    \]
    with $Z_{d,\alpha}$ chosen such that $\int\mu_\alpha = 1$ is a probability distribution such that $L_\alpha\mu_\alpha = 0$. Observe that away from the origin, the Fourier transform of $\mu_\alpha$ is smooth and rapidly decaying at infinity. The singularity at~$\xi=0$ behaves like $|\xi|^\alpha$ at principal order which yields that $\mu_\alpha(v)$ should decay as $|v|^{-\alpha-d}$ when~$|v|\to\infty$. Actually, one has the following more precise bounds { coming from \cite[Theorem 3.1]{bogdan_2003_harnack} and references in the proof} (see also the references in \cite{aceves_2019_fractional}). There are positive constants $C_1 = C_1(\alpha,d)>0$ and $C_2 = C_2(\alpha,d)>0$ such that for all $v\in\RR^d$ one has
    \begin{equation}\label{eq:mualpha_est1}
        C_1^{-1}\ \leq\ (|v|^{d+\alpha}+1)\mu_\alpha(v)\ \leq\ C_1\,,
    \end{equation}
    and
    \begin{equation}\label{eq:mualpha_est2}
        C_2^{-1}|v|\,\ \leq\ (|v|^{2+d+\alpha}+1)|\nabla_v\mu_\alpha(v)|\ \leq\ C_2|v|\,.
    \end{equation}
    
    In the following, given some measurable non-negative function $\nu\equiv\nu(v)$ we denote by $L^2_{v}(\nu)$  and~$L^2_{x,v}(\nu)$ the spaces of measurable functions $g$ of respectively the $v$ and the $(x,v)$ variables such that $|g|^2\nu$ is integrable. We endow these spaces with their canonical scalar product and norm. We also introduce the corresponding Sobolev space $H^1_{x,v}(\nu)$ associated with the norm
    $$
    \|g\|_{H^1_{x,v}(\nu)}^2\,=\,\|g\|_{L^2_{x,v}(\nu)}^2 + \|\nabla_x g\|_{L^2_{x,v}(\nu)}^2 + \|\nabla_v g\|_{L^2_{x,v}(\nu)}^2\, .
    $$ 
    Finally given an integrable function $g$, we denote $\lla g\rra = \iint_{\TT^d\times\RR^d} g(x,v) \, \dd v \, \dd x$ the global mass of $g$.
    The main result of this paper is the following.
    
    \medskip
    \begin{theorem}\label{t:main}
        Let $f$ solve the kinetic L\'evy-Fokker-Planck equation \eqref{eq:kinfracFP} with initial data~$f^\text{in}\in H^1_{x,v}(\mu_\alpha^{-1})$.  Then, for all $t\geq0$ one has
        \[
        \|f(t) - \lla f^\text{in}\rra\,\mu_\alpha\|_{H^1_{x,v}(\mu_\alpha^{-1})}\ \leq\ C\,\|f^\text{in} - \lla f^\text{in}\rra\,\mu_\alpha\|_{H^1_{x,v}(\mu_\alpha^{-1})}\,e^{-\lambda t}\,
        \]
        for some constant $C\geq1$ and $\lambda>0$ depending only on $d$ and $\alpha$. 
    \end{theorem}

    \medskip
    Let us mention that these results have been obtained as a preliminary step towards the conception and analysis of numerical schemes preserving the long-time behavior of these equations. This topic is an ongoing work \cite{ayi_2020_fractional} in the spirit of what has previously been done in~\cite{dujardin_2018_coercivity,bessemoulin_2018_hypocoercivity} in the case of the classical Fokker-Planck equation. The compatibility of our schemes with anomalous diffusion limit will also be investigated (see \cite{crouseilles_2016_numerical} for more details).

        Before going into the analysis of our problem, let us recall that results on large time behavior of solutions to the homogeneous version of~\eqref{eq:kinfracFP}, namely $\partial_tf(t,v) = L_\alpha f(t,v)$, have been obtained in~\cite{gentil_2008_levy} in spaces of type $L^2_v(\mu_\alpha^{-1})$ (among others) and later in~\cite{tristani_2015_FFP} in larger Lebesgue spaces. Notice that the presence of the transport operator in our equation~\eqref{eq:kinfracFP} makes the analysis more intricate and requires the use of hypocoercivity techniques. In the present note, we use $H^1$ type hypocercivity as presented in~\cite{villani_2009_hypocoercivity} or \cite{herau_2018_introduction} for example. {Note also that fractional hypocoercivity has already been studied recently in~\cite{bouin_2019_fractional} where a $L^2$-hypocoercivity approach is developed. In this sense, their framework is quite different, note however that it is also more general than ours (in terms of phase space and linear operators). Their results in particular imply an exponential convergence towards equilibrium in the torus in $L^2$ for our models.}
        
        In the same spirit of our work, let us also mention the paper~\cite{herau_2018_frac} in which some hypoelliptic estimates are obtained on the non homogeneous fractional Kolmogorov equation (there is no drift term in the studied equation). The method of proof is quite close (based on the use  of weighted Lyapunov functional) but the final goal is different in the latter since the main concern is about regularization properties of the equation and not convergence towards the equilibrium.

        In the present study we focus on a good understanding of the structure of the L\'evy-Fokker-Planck operator since we endeavour to carry out our computations as simply as possible in order to adapt our analysis to a discrete framework in  \cite{ayi_2020_fractional}. In particular let us point out  that we do not need  fractional derivatives in our Lyapunov functionals and our proof does not rely on Fourier transform. In this sense our method differs completely from that of \cite{herau_2018_frac} and the recent \cite{bouin_2019_fractional} in which a mode by mode analysis is developed.
    
    \smallskip
    \noindent{\it Outline of the note.}
    From Section~\ref{s:bilinear} to Section~\ref{s:interpolation}, we carry out the analysis of the  properties of the L\'evy-Fokker-Planck operator that will be useful for proving our main result. Then, the proof of Theorem~\ref{t:main} is done in Section~\ref{s:proof}. In the last section we state and prove the equivalent of Theorem~\ref{t:main} for the BGK equation with heavy-tailed equilibrium.  
    
    \smallskip
    
    \noindent{\it Notations.} For simplicity, in the subsequent proofs, we denote by $C$ a positive constant depending only on fixed numbers (including $d$ and $\alpha$) and its value may change from line to line.
    
    \section{The L\'evy-Fokker-Planck operator as bilinear form}\label{s:bilinear}
    
    The following quite simple decomposition is actually one of the key elements of our hypocoercive analysis carried out in Section~\ref{s:proof}. Compared to the non-fractional case, we here have a lack of symmetry of our operator in $L^2_v(\mu_\alpha^{-1})$ and the following splitting is very helpful to simplify the computations. Moreover, in the non-fractional case, there is a gain of weight in velocity which comes from the particular form of the gradient of the Gaussian equilibrium. Even though we no longer have such a gain in our case, we are still able to close our estimates thanks to the following splitting.
    
    \smallskip
    \begin{e-proposition} \label{prop:decomp}
        One has the decomposition
        \begin{equation*}
            -\lla L_\alpha f,\,g\rra_{L^2_v(\mu_\alpha^{-1})}\ =\ \cS_v(f,g)\,+\,\cA_v(f,g)\,,
        \end{equation*}
        where $\cS_v$ and $\cA_v$ are bilinear forms that are respectively symmetric and skew-symmetric and defined by
        \begin{equation*}
            \cS_v(f,g)\ =\ \frac{C_{d,\alpha}}{2}\iint_{\RR^d\times\RR^d}\frac{\left[({f}\mu_\alpha^{-1})(v)-({f}\mu_\alpha^{-1})(w)\right]\left[({g}\mu_\alpha^{-1})(v)-({g}\mu_\alpha^{-1})(w)\right]}{|v-w|^{d+\alpha}}\mu_\alpha(v)\, \dd w\,\dd v\,,
        \end{equation*}
        and
        \begin{multline*}
            \cA_v(f,g)\ =\ \frac{C_{d,\alpha}}{2}\iint_{\RR^d\times\RR^d}\frac{({f}\mu_\alpha^{-1})(w)({g}\mu_\alpha^{-1})(v)-({f}\mu_\alpha^{-1})(v)({g}\mu_\alpha^{-1})(w)}{|v-w|^{d+\alpha}}\mu_\alpha(v)\, \dd w\,\dd v\\\,+\, \frac{1}{2}\int_{\RR^d}(f\,v\cdot\nabla_v({g}\mu_\alpha^{-1})\,-\,g\,v\cdot\nabla_v({f}\mu_\alpha^{-1}))\, \dd v\,,
        \end{multline*}
         where $C_{d,\alpha}$ is defined in~\eqref{eq:fraclap_int1}. 
    \end{e-proposition}
    We skip the proof of this proposition since it is based on simple computations using the formula~\eqref{eq:fraclap_int1}, integration by parts and the fact that $L_\alpha \mu_\alpha=0$.

    Observe that a direct consequence of the Cauchy-Schwarz inequality is 
    \begin{equation}\label{eq:cauchy_schwarz}
        \cS_v(f,g)\ \leq\ \cS_v(f,f)^{1/2}\,\cS_v(g,g)^{1/2}\,,
    \end{equation}
    for $f,g\in D(L_\alpha)$. Moreover, the symmetric form $\cS_v$ is non-negative and $\cS_v(f,f)$ vanishes when $f\mu_\alpha^{-1}$ is constant. This yields that the nullspace of $L_\alpha$ is exactly given by $\RR\mu_\alpha$.  From there the orthogonal projection $\Pi$ onto the nullspace of $L_\alpha$ is given by 
    \begin{equation*}
        (\Pi g)(v)\ =\ \left(\int_{\RR^d}g(w)\,\dd w\right)\,\mu_\alpha(v)\,.
    \end{equation*}

    \section{Coercivity results for the L\'evy-Fokker-Planck operator}\label{s:coercive}

    One has the following coercivity result taken from \cite[Theorem 2]{gentil_2008_levy} and originating from \cite{chafai_2004_entropies}. { While the previous references derive the inequality via a semigroup approach, let us mention that an elementary analytical proof is given by Wang \cite{wang_2014_simple} and came to our attention thanks to \cite{bouin_2019_fractional}.}
    
    \begin{lemma}[\cite{chafai_2004_entropies}, \cite{gentil_2008_levy}, \cite{wang_2014_simple}...]\label{eq:estim}
        There is a constant $C_P\equiv C_P(\alpha,d)>0$ such that for all $f\in D(L_\alpha)$, 
        \begin{equation} \label{eq:coerciv}
            \|f-\Pi f\|_{L^2_v(\mu_\alpha^{-1})}^2\ \leq\ C_P\,\cS_v(f,f) \, .
        \end{equation}
    \end{lemma}
    
    We now show that the dissipation $\cS_v(f,f)$  also provides some fractional Sobolev regularity. We introduce the fractional Sobolev space $H^{s}_v$ with $s\in(0,1)$ with norm defined by $\|\cdot\|_{H^{s}_v}^2 = \|\cdot\|_{L^2_v}^2\,+\,\|\cdot\|_{\dot{H}^{s}_v}^2$ where the homogeneous Sobolev norm is given by
    $
    \|g\|_{\dot{H}^{s}_v}^2 :=\ \|(-\Delta)^{s/2}g\|_{L^2_v}^2
    $.
    One can prove that there exists a positive constant $\widetilde{C}_{d,s}$ such that 
        \begin{equation} \label{eq:fracSob}
            \|g\|_{\dot{H}^{s}_v}^2 = \widetilde{C}_{d,s} \iint_{\RR^d\times\RR^d} {|f(v)-f(w)|^2}{|v-w|^{-(d+2s)}} \, \dd w \, \dd v \,.
        \end{equation}
    \smallskip
    \begin{lemma} \label{lem:estimDalpha}
        There exists $C_R\equiv C_R(\alpha,d)>0$ such that for all $f\in D(L_\alpha)$,
        \[
        \cS_v(f,f)\ \geq\  C_R^{-1}\,\left(\|  f\mu_\alpha^{-1/2} \|_{\dot H^{\alpha/2}_v}^{2} -  \|  f \mu_\alpha^{-1/2} \|_{L^2_v}^{2}\right).
        \]
    \end{lemma}
    \noindent\emph{Proof. }
    Using that $(a+b)^2\geq a^2/2-b^2$, we have
    $$
    \cS_v(f,f) 
    \ge \frac{C_{d,\alpha}}{2} \iint_{|v-w| \le 1}   \frac{|(\mu_\alpha^{-1} f)(v) - (\mu_\alpha^{-1} f)(w)|^2}{|v-w|^{d+\alpha}} \, \mu_\alpha(v) \, \dd w \, \dd v
    \geq  \frac{C_{d,\alpha}}{2} \bigg(\frac12I_1 - I_2\bigg)\,.
    $$
    The first term is 
    \[
    \begin{array}{rcl}
        I_1 
        \, &=&\,  \ds\iint_{|v-w| \le 1} \frac{|(\mu_\alpha^{-1/2} f) (v) - (\mu_\alpha^{-1/2} f)(w)|^2}{|v-w|^{d+\alpha}} \, \dd w \, \dd v \\
        &=& \ds{ \widetilde{C}_{d,{\frac{\alpha}{2}}}^{-1}} \| \mu_\alpha^{-1/2} f \|_{\dot H^{\alpha/2}_v}^{2} - \iint_{|v-w| \ge 1 } \frac{|(\mu_\alpha^{-1/2} f) (v) - (\mu_\alpha^{-1/2} f)(w)|^2}{|v-w|^{d+\alpha}} \, \dd w \, \dd v\\
        &\geq& \ds{ \widetilde{C}_{d,{\frac{\alpha}{2}}}^{-1}}\| \mu_\alpha^{-1/2} f \|_{\dot H^{\alpha/2}_v}^{2}
        - {C} \| \mu_\alpha^{-1/2} f \|_{L^2_v}^{2}
    \end{array}
    \]
    where { $\widetilde{C}_{d,{\frac{\alpha}{2}}}$ is defined in~\eqref{eq:fracSob}} and for the last inequality, we used the integrability of $|v-w|^{-d-\alpha} { \mathds{1}_{|v-w| \ge 1}}$, once in $v$ and once in $w$. The second term is
    \[
    I_2\ =\  \iint_{|v-w| \le 1} \frac{ |\mu_\alpha^{1/2} (v) - \mu_\alpha^{1/2} (w)|^2}{|v-w|^{d+\alpha}} \,  |f(w)|^2 |\mu_\alpha^{-1} (w)|^2 \, \dd w \,  \dd v \, .
    \]
    To treat $I_2$, we use Taylor formula to write
    \[
    I_2 = \iint_{|w| \le 1} \frac{\left|\int_0^1 \nabla(\mu_\alpha^{1/2})(v-\theta w) \cdot w \, \dd\theta\right|^2}{|w|^{d+\alpha}} \mu_\alpha^{-1}(v-w) |f(v-w) \mu_\alpha^{-1/2}(v-w)|^2 \, \dd w \, \dd v\,. 
    \]
    Performing now the changes of variables $v \to v- \theta w$ and then $ \theta \to 1-\theta$, we get:
    \[
    I_2 \le \iint_{|w| \le 1} \int_0^1\frac{|\nabla(\mu_\alpha^{1/2})(v)|^2}{|w|^{d+\alpha-2}} \mu_\alpha^{-1}(v-\theta w) |f(v-\theta w)\mu_\alpha^{-1/2}(v-\theta w)|^2 \, \dd\theta \,  \dd w \,  \dd v\,.
    \]
    Notice that, using~{\eqref{eq:mualpha_est1}} and since $|w| \le 1$, we have 
    $
    \mu_\alpha^{-1}(v-\theta w) \leq {C}(1 + |v|^{d+\alpha}).
    $ 
    Then, using~\eqref{eq:mualpha_est2}, one can prove that 
    $
    |\nabla(\mu_\alpha^{1/2})(v)|^2\mu_\alpha^{-1}(v-\theta w) \leq {C}. 
    $ 
    Consequently, we obtain 
    $$
    I_2 \le {C} \iint_{|w| \le 1} \int_0^1 {\frac{1}{|w|^{d+\alpha-2}}} |f(v-\theta w)\mu_\alpha^{-1/2}(v-\theta w)|^2 \, \dd\theta \, \dd w \, \dd v
    $$
    and thus performing a change of variable 
    $
    I_2 \le {C} \| f  \mu_\alpha^{-1/2} \|_{L^2_v}^{2}. 
    $ This ends the proof.\qed
    
    \begin{e-proposition} \label{prop:CF}
        There is $C_F\equiv C_F(\alpha,d)$ such that for all $f\in D(L_\alpha)$,
        \begin{equation} \label{eq:coerciv2}
            \|(f- \Pi f)\mu_\alpha^{-1/2}\|^2_{H^{\alpha/2}_v}\ \leq\ C_F\,\cS_v(f,f)\,.
        \end{equation}
    \end{e-proposition}
    \noindent\emph{Proof. }
    Let us now summarize the estimates that we have obtained in the two previous lemma. We have~$\cS_v(f,f)\geq C_P^{-1} \|f - \Pi f\|^2_{L^2_v(\mu_\alpha^{-1/2})}$ 
    and $\cS_v(f,f)\ \geq\  {C_R^{-1}}\,\big(\|  f\mu_\alpha^{-1/2} \|_{\dot H^{\alpha/2}_v}^{2} -  \|  f \mu_\alpha^{-1/2} \|_{L^2_v}^{2}\big)$. Moreover, one can notice that $
    \cS_v(f,f)\ =\ \cS_v(f-\Pi f,f-\Pi f)$. As a consequence, an appropriate convex combination of the two previous inequalities shows \eqref{eq:coerciv2}.
    \qed
    
    \section{An interpolation inequality}\label{s:interpolation}
    
    In this section we prove an interpolation result which is crucial in the proof of Theorem~\ref{t:main}.
    \begin{e-proposition}
        For all $\ve>0$, there is $K(\ve)\equiv K(\ve,\alpha,d)>0$ such that 
        \begin{equation}\label{eq:interp}
            \|\nabla_vf\|_{L^2_v(\mu_\alpha^{-1})}^2\ \leq\ K(\ve) \left(\cS_v(f,f)\,+\,\|\Pi f\|_{L^2_v(\mu_\alpha^{-1})}^2\right)+\, \ve\, {C_F}\,\cS_v(\nabla_vf,\nabla_vf)
        \end{equation}
        where the constant $C_F$ is defined in Proposition~\ref{prop:CF}. 
    \end{e-proposition}
    
    \noindent\emph{Proof. }
    One can use the chain rule and an interpolation of $\dot H^1_v$ between $\dot H^{\alpha/2}_v$ and $\dot H^{1+\alpha/2}_v$ (easily shown in Fourier variables) to get
    \[
    \begin{array}{rcl}
        \|(\nabla_v f)\mu_\alpha^{-1/2}\|^2_{L^2_{v}} &\leq&{2}\| \nabla_v (f \mu_\alpha^{-1/2})\|^2_{L^2_{v}} +{ 2}\|f (\nabla_v \mu_\alpha^{-1/2})\|^2_{L^2_{v}} \\[.5em]
        &\leq& K(\eps)\, \| f \mu_\alpha^{-1/2}\|^2_{\dot H^{\alpha/2}_v} + \eps\, \|\nabla_v(f \mu_\alpha^{-1/2})\|^2_{\dot H^{\alpha/2}_v} + {2}\|f (\nabla_v \mu_\alpha^{-1/2})\|^2_{L^2_{v}}\\[.5em]
        &\leq& K(\eps)\, \|f \mu_\alpha^{-1/2}\|^2_{H^{\alpha/2}_v} + \eps\, \|(\nabla_v f)\mu_\alpha^{-1/2}\|^2_{H^{\alpha/2}_v}+ C\,\|f \mu_\alpha^{-1/2}\|^2_{L^2_{v}}\\[.5em]
        &\leq& K(\eps)\, \|f \mu_\alpha^{-1/2}\|^2_{H^{\alpha/2}_v} + \eps\, \|(\nabla_v f)\mu_\alpha^{-1/2}\|^2_{H^{\alpha/2}_v}
    \end{array}
    \]
    {up to changing the value of $K(\eps)$} and where we used the fact that $|(\nabla_v \mu_\alpha) \mu_\alpha^{-1}| \in L^\infty(\RR^d)$ to bound the third term.
    Now observe that
    \[
    \|f \mu_\alpha^{-1/2}\|^2_{H^{\alpha/2}_v}\ \leq\ 2\bigg(\|(f-\Pi f) \mu_\alpha^{-1/2}\|^2_{H^{\alpha/2}_v}+\|(\Pi f) \mu_\alpha^{-1/2}\|^2_{H^{\alpha/2}_v}\bigg)\,,
    \]
    and that 
    $
    \|(\Pi f)  \mu_\alpha^{-1/2}\|_{H^{\alpha/2}_v}\ =\ \|\Pi f\|_{L^2_v(\mu_\alpha^{-1})}\|\mu_\alpha^{1/2}\|_{H^{\alpha/2}_v}
    $
    {with $\|\mu_\alpha^{1/2}\|_{H^{\alpha/2}_v} \le C$ since $\mu_\alpha^{1/2} \in H^1_v$ from~\eqref{eq:mualpha_est1} and~\eqref{eq:mualpha_est2}}. 
    Moreover, one has $\nabla_v f = \nabla_v f - \Pi\nabla_v f$. One can conclude by using \eqref{eq:coerciv2} twice.
    \qed

    \section{Proof of Theorem~\ref{t:main}}\label{s:proof}
    
    Up to changing $f^\text{in}$ by {$f^\text{in} - \lla f^\text{in} \rra \mu_\alpha$},  we assume that  $\iint_{\TT^d\times\RR^d}f(t,x,v) \, \dd v \, \dd x \ =\ 0$ at $t=0$, so that by conservation it also holds for all time $t>0$. We introduce a new norm on the weighted Sobolev space~$H^1_{x,v}(\mu_\alpha^{-1})$. It is defined by
    \begin{equation} \label{def:triplenorm}
        \Nt f\Nt^2 
        \ =\  \|f\|^2_{L^2_{x,v}(\mu_\alpha^{-1})} \,+\, a\|\nabla_x f \|^2_{L^2_{x,v}(\mu_\alpha^{-1})} \,+\, b\, \|\nabla_v f \|^2_{L^2_{x,v}(\mu_\alpha^{-1})}  + 2\,c\, \langle \nabla_x f , \nabla_v f \rangle_{L^2_{x,v}(\mu_\alpha^{-1})}\,,
    \end{equation}
    where $a$, $b$ and $c$ are positive constants to be determined later on. Observe that as soon as {$c^2<ab$}, one has that $\Nt\cdot\Nt$ is equivalent to $\|\cdot\|_{H^1_{x,v}(\mu_\alpha^{-1})}$. Let us note that the commutators $[\nabla_x, v\cdot\nabla_x]$ and $[\nabla_x, L_\alpha]$ vanish while $[\nabla_v, v\cdot\nabla_x]\ =\ \nabla_x$ and $[\nabla_v, L_\alpha]\ =\ \nabla_v$. Also observe that $v\cdot\nabla_x$ is skew-symmetric in~$L^2_{x,v}(\mu_\alpha^{-1})$. {Let us estimate the evolution of each term appearing in the new norm defined in \eqref{def:triplenorm} for~$f$ a solution of~\eqref{eq:kinfracFP} with initial data $f^{in}$ satisfying $\lla f^{in} \rra$ = 0}. In the following the notation $\cS_{x,v}$ denotes the integral of {$\cS_{v}$} in the $x$ variable. One has
    \begin{align*}
        \ds\frac12\,\frac{\dd}{\dd t}\|f\|^2_{L^2_{x,v}(\mu_\alpha^{-1})} &= -\cS_{x,v}(f,f)\, ,\\
        \ds\frac12\,\frac{\dd}{\dd t}\|\nabla_xf\|^2_{L^2_{x,v}(\mu_\alpha^{-1})} &= -\cS_{x,v}(\nabla_x f,\nabla_x f)\, ,\\
        \ds\frac12\,\frac{\dd}{\dd t}\|\nabla_vf\|^2_{L^2_{x,v}(\mu_\alpha^{-1})} &= -\cS_{x,v}(\nabla_v f,\nabla_v f)+\,  \|\nabla_vf\|_{L^2_{x,v}(\mu_\alpha^{-1})}^2\,-\,\lla\nabla_xf,\nabla_vf \rra_{L^2_{x,v}(\mu_\alpha^{-1})}\, ,\\
        \ds\,\frac{\dd}{\dd t}\lla\nabla_xf,\nabla_vf \rra_{L^2_{x,v}(\mu_\alpha^{-1})} &= -\|\nabla_xf\|_{L^2_{x,v}(\mu_\alpha^{-1})}^2-2\,\cS_{x,v}(\nabla_x f,\nabla_v f) + \lla\nabla_xf,\nabla_vf \rra_{L^2_{x,v}(\mu_\alpha^{-1})} \, .
    \end{align*} 
    Notice here that the keystone of the proof of the last equality is the splitting obtained in Proposition~\ref{prop:decomp} { and the Hilbertian setting. Indeed given any $g =e^{tL_\alpha}g_0$ and operators $A$ and $B$, one has formally that $\frac{\dd}{\dd t}\lla Ag,Bg\rra = \lla [A,L_\alpha]g,Bg\rra + \lla Ag,[B,L_\alpha]g\rra -2\cS_{x,v}(Ag,Bg)$. Therefore the skew symmetric part of the operator $L_\alpha$ only appears in commutators. This observation enables us to avoid loss of moments in velocities in forthcoming estimates which one would face with bad rearrangements of the terms.}
    By gathering all the previous estimates one gets 
    \begin{multline*}\label{}
        \frac 12\frac{\dd}{\dd t}\Nt f\Nt^2\ =\ -\cS_{x,v}(f,f) - a\cS_{x,v}(\nabla_x f,\nabla_x f)-b\cS_{x,v}(\nabla_v f,\nabla_v f) -c\|\nabla_xf\|_{L^2_{x,v}(\mu_\alpha^{-1})}^2\\[.75em]
        +\,  b\|\nabla_vf\|_{L^2_{x,v}(\mu_\alpha^{-1})}^2\,-\,b\lla\nabla_xf,\nabla_vf \rra_{L^2_{x,v}(\mu_\alpha^{-1})}\\[.75em]
        -2c\,\cS_{x,v}(\nabla_x f,\nabla_v f) + c\lla\nabla_xf,\nabla_vf \rra_{L^2_{x,v}(\mu_\alpha^{-1})} \, .
    \end{multline*}
    The first four terms are dissipation terms and the last  four terms are remainder terms. Let us control the latter by the former ones. By integrating \eqref{eq:cauchy_schwarz} in $x$ and using Young's inequality one gets 
    \[
    \left|2c\,\cS_{x,v}(\nabla_x f,\nabla_v f)\right|\ \leq\  \frac{2\,c^2}{b}\,\cS_{x,v}(\nabla_x f,\nabla_x f) + \frac{b}{2}\,\cS_{x,v}(\nabla_v f,\nabla_v f)\,.
    \]
    Then since $\int \nabla_vf \dd v = 0$, one has
    \[
    \begin{array}{rcl}
        b\left|\lla\nabla_xf,\nabla_vf \rra_{L^2_{x,v}(\mu_\alpha^{-1})}\right|& =& b\left|\lla\nabla_xf-\Pi(\nabla_xf),\nabla_vf \rra_{L^2_{x,v}(\mu_\alpha^{-1})}\right|\\[.75em]
        & \leq& b\, {C_P^{1/2}}\cS_{x,v}(\nabla_x f,\nabla_x f)^{1/2}\|\nabla_vf \|_{L^2_{x,v}(\mu_\alpha^{-1})}\\[.75em]
        & \leq& \ds\frac{b\,{C_P}}{2}\,\cS_{x,v}(\nabla_x f,\nabla_x f) + \frac{b}{2}\,\|\nabla_vf \|_{L^2_{x,v}(\mu_\alpha^{-1})}^2\,,
    \end{array}
    \]
    where we used \eqref{eq:coerciv}. Similarly
    \[
    c\left|\lla\nabla_xf,\nabla_vf \rra_{L^2_{x,v}(\mu_\alpha^{-1})}\right|\ \leq\ \ds\frac{c^2 {C_P}}{2b}\,\cS_{x,v}(\nabla_x f,\nabla_x f) + \frac{b}{2}\,\|\nabla_vf \|_{L^2_{x,v}(\mu_\alpha^{-1})}^2\,.
    \] 
    For the last remainder term we use \eqref{eq:interp} integrated in $x$, namely
    \[
    \|\nabla_vf\|_{L^2_{x,v}(\mu_\alpha^{-1})}^2\ \leq\ K(\ve) \left(\cS_{x,v}(f,f)\,+\,\|\Pi f\|_{L^2_{x,v}(\mu_\alpha^{-1})}^2\right)+\, \ve\,{C_F}\, \cS_{x,v}(\nabla_vf,\nabla_vf)\,.
    \]
    We can use the Poincar\'e inequality on the torus (since $f$ is mean-free) and the Jensen inequality to get~
    $
    \|\Pi f\|_{L^2_{x,v}(\mu_\alpha^{-1})}^2\ \leq\ {\widetilde{C}_P}\|\nabla_xf\|_{L^2_{x,v}(\mu_\alpha^{-1})}^2
    $ 
    {where $\widetilde{C}_P \equiv \widetilde{C}_P(d)$ is the Poincar\'e constant of the $d$-dimensional torus.}
    Thus eventually, one has 
    \begin{equation}\label{e:entropdissip}
        {\frac12} \frac{\dd}{\dd t}\Nt f\Nt^2 + D(f,f)\leq 0\,,
    \end{equation}
    where the dissipation is given by
    \begin{multline*}
        D(f,f)\ =\ (1 - 2\,b\,K(\eps))\cS_{x,v}(f,f) \,+\,{\left(a-\frac{c^2}{b}\left(2+\frac{C_P}{2}\right)-\frac{b C_P}{2}\right)}\cS_{x,v}(\nabla_x f,\nabla_x f)\\
        +\left({\frac{b}{2}}-2b\eps{C_F}\right)\,\cS_{x,v}(\nabla_v f,\nabla_v f) +\left(c-2b{\widetilde{C}_P}K(\eps)\right)\,\|\nabla_xf\|_{L^2_{x,v}(\mu_\alpha^{-1})}^2\,.
    \end{multline*}
    
    Now choose consecutively $\eps,b,c$ and $a$ such that {$0<\eps<1/(4C_F)$}, $0<b<1/(2K(\eps))$, $c>2b{\widetilde{C}_P}K(\eps)$ and finally $a$ large enough so that $a>{c^2}\left(2+{C_P}/{2}\right)/b+{b C_P}/{2}$. It yields that the dissipation is non-negative and even that there is a constant $\lambda>0$ (depending on $a,b,c,\eps$) such that
    $
    D(f,f)\geq\lambda \Nt f\Nt^2
    $. 
    By a Gronwall type argument we have that $\Nt f(t)\Nt$ decays exponentially to $0$ when $t\to\infty$. \qed

    \section{The case of the heavy-tailed BGK equation}\label{s:bgk}
    In this last section we consider another simple kinetic model
    \begin{equation}\label{eq:kinfracBGK}
        \partial_t f + v\cdot\nabla_x f\ =\ \Pi_Mf - f\,,\qquad\text{with}\quad (\Pi_Mf)(t,x,v)\ =\ M(v)\,\int_{\RR^d}f(t,x,w)\,\dd w\,,
    \end{equation}
    for which the local equilibrium satisfies the following assumptions
    \begin{equation}\label{eq:hypotheses}
        M(v)>0\,, \qquad \int_{\RR^d}M = 1\,,\quad\text{and}\quad \nabla_v\ln(M)\in L^\infty\,.
    \end{equation}
    This allows for heavy-tailed distributions, namely $M$ such that $M(v) \sim_{|v|\to\infty}|v|^{-d-\alpha}$ with $\alpha\in(0,2)$.
    \begin{theorem}\label{t:bgk}
        Assume that \eqref{eq:hypotheses} holds and let $f$ solve the BGK equation \eqref{eq:kinfracBGK} starting from the initial data $f^\text{in}\in H^1_{x,v}(M^{-1})$. Then, for all $t\geq0$ one has
        \begin{equation*}
            \|f(t) - \lla f^\text{in}\rra\,M\|_{H^1_{x,v}(M^{-1})}\ \leq\ {C}\,\|f^\text{in} - \lla f^\text{in}\rra\,M\|_{H^1_{x,v}(M^{-1})}\,e^{{-\lambda} t}\,
        \end{equation*}
        for some constant {$C\geq1$} and {$\lambda>0$} depending only {on $d$ and $\|\nabla_v\ln(M)\|_{L^\infty}$}. 
    \end{theorem}
    
    The proof is similar and simpler than that of Theorem~\ref{t:main}. We skip many details as the reader may go back to the proof of Theorem~\ref{t:main} in order to recover them.
    
    \noindent\emph{Proof of Theorem~\ref{t:bgk}. } {Consider $f$ a solution to~\eqref{eq:kinfracBGK} with initial data $f^{in}$ satisfying~$\lla f^{in} \rra=0$.} Let us observe that the commutators $[\nabla_x, v\cdot\nabla_x]$ and $[\nabla_x, \Pi_M]$ vanish while $[\nabla_v, v\cdot\nabla_x]\ =\ \nabla_x$ and also~$[\nabla_v, \Pi_M]\ =\ \nabla_v\ln(M)\Pi_M$. Now with this in mind, and defining the triple norm of $f$ as in~\eqref{def:triplenorm} with~$\mu_\alpha$ replaced by $M$, one gets  
    \begin{multline*}\label{}
        \frac 12\frac{\dd}{\dd t}\Nt f\Nt^2\ =\ -\|f-\Pi_Mf\|^2_{L^2_{x,v}(M^{-1})} - a\|\nabla_xf-\Pi_M\nabla_xf\|^2_{L^2_{x,v}(M^{-1})}-b\|\nabla_vf\|^2_{L^2_{x,v}(M^{-1})} -c\|\nabla_xf\|_{L^2_{x,v}(M^{-1})}\\[.75em]
        +\,  b\lla\nabla_v\ln(M)\Pi_Mf,\nabla_vf\rra_{L^2_{x,v}(M^{-1})}\,-\,b\lla\nabla_xf,\nabla_vf \rra_{L^2_{x,v}(M^{-1})}\\[.75em]
        -2c\lla\nabla_xf,\nabla_vf\rra_{L^2_{x,v}(M^{-1})} + c\lla\nabla_v\ln(M)\Pi_Mf,\nabla_xf\rra_{L^2_{x,v}(M^{-1})}\,.
    \end{multline*}
    {First, we notice that $\lla\nabla_xf,\nabla_vf \rra_{L^2_{x,v}(M^{-1})}=\lla\nabla_xf-\Pi_M \nabla_xf,\nabla_vf \rra_{L^2_{x,v}(M^{-1})}$ to deal with the third and fourth remainder terms with Cauchy-Schwarz inequality.} The last remainder term requires some special care. Indeed, observe that since~$\lla\nabla_v\ln(M)\Pi_Mf,\Pi_Mg\rra$ vanishes for any $g$, one thus has
    \[
    \lla\nabla_v\ln(M)\Pi_Mf,\nabla_xf\rra_{L^2_{x,v}(M^{-1})}\ \leq\ \|\nabla_v\ln(M)\|_{L^\infty}\,\|\Pi_Mf\|_{L^2_{x,v}(M^{-1})}\,\|\nabla_xf - \Pi_M\nabla_xf\|_{L^2_{x,v}(M^{-1})} \, .
    \]
    {We also have that 
        \[
        \lla\nabla_v\ln(M)\Pi_Mf,\nabla_vf\rra_{L^2_{x,v}(M^{-1})} \ \leq \  \|\nabla_v\ln(M)\|_{L^\infty}\,\|\Pi_Mf\|_{L^2_{x,v}(M^{-1})}\,\|\nabla_vf \|_{L^2_{x,v}(M^{-1})} \, .
        \]}    
    Finally, we recall that $\|\Pi_Mf\|_{L^2_{x,v}(M^{-1})}\leq {\widetilde{C}_P}\,\|\nabla_xf\|_{L^2_{x,v}(M^{-1})}$ with $ {\widetilde{C}_P}$ the Poincar\'e constant of the $d$-dimensional torus. Then using four times Young's inequality with well chosen weights, one obtains \eqref{e:entropdissip} with the dissipation 
    \begin{multline*}
        D(f,f)\ =\ \|f-\Pi_Mf\|^2_{L^2_{x,v}(M^{-1})}+ \left(a-b-4c^2/b-C_{M}c/2\right)\|\nabla_xf-\Pi_M\nabla_xf\|^2_{L^2_{x,v}(M^{-1})}\\+\left(b-b/4-b/4-b/4\right)\|\nabla_vf\|^2_{L^2_{x,v}(M^{-1})} +\left(c - bC_{M}-c/2\right)\|\nabla_xf\|_{L^2_{x,v}(M^{-1})}
    \end{multline*}
    with {$C_{M} = \|\nabla_v\ln(M)\|_{L^\infty}^2 {\widetilde{C}_P}^2$}. One concludes as in Theorem~\ref{t:main} after choosing any $b>0$, $c>2\,b\,C_{M}$ and finally $a>4c^2/b+b+C_{d,M}c/2$.
    \qed

    
    
    
    \section*{Acknowledgements}
    Maxime Herda thanks the LabEx CEMPI (ANR-11-LABX-0007-01). H\'el\`ene Hivert thanks the European Research Council (ERC) under the European Union's Horizon $2020$ research and innovation programme (ERC starting grant MESOPROBIO n${}^\text{o}$ $639638$). Isabelle Tristani thanks the ANR EFI:  ANR-17-CE40-0030  and the ANR SALVE: ANR-19-CE40-0004 for their support.
    
    This work is part of a collaborative research project that was initiated for the Junior Trimester Program in Kinetic Theory at the Hausdorff Research Institute for Mathematics in Bonn. Part of the work was carried out during the time spent at the institute. The authors are grateful for this opportunity and warmly acknowledge the HIM for the financial support and the hospitality they benefited during their stay.

\end{document}